\documentclass[11pt,a4paper,openany,twoside]{article}
\usepackage{amsmath}
\usepackage{amssymb}
\usepackage[amsmath,thmmarks]{ntheorem}
\usepackage{fancyhdr}
\usepackage[top=2.85cm,bottom=2.8cm,left=37mm,right=37mm,includehead,includefoot]{geometry}
\theoremseparator{.\hspace{0.2em}}

\newtheorem{lemma}{\bf Lemma}

\newtheorem{example}{\bf Example}
\newtheorem{claim}{\bf Claim}
\newtheorem{maintheorem}{\bf Main Theorem}

\theoremstyle{nonumberplain}
\theorembodyfont{} \theoremindent0em

\theorembodyfont{\normalfont \rm } \theoremindent0em
\theoremseparator{\hspace{1em}}
\def\simlin{\sim_{\text{\rm lin}}}

\def\simnum{\sim_{\text{\rm num}}}
\theoremsymbol{$\blacksquare$}
\newenvironment{proof}[1][Proof.]{\begin{trivlist}
\item[\hskip \labelsep {\bfseries #1}]}{\end{trivlist}}
\def\roundup#1{\ulcorner{#1}\urcorner}

\def\Co#1{{\mathcal{O}}_{#1}}
 \allowdisplaybreaks[4]
\def\baselinestretch{1.3}

\begin{document}
\title{The sharp lower bound for the volume of 3-folds of general type with $\chi(\Co{X})=1$}
\author{Lei Zhu\footnotemark[1]\\
\thanks{Supported by Fudan Graduate Students' Innovation Projects (EYH5928004).}
{\footnotesize \it School of Mathematical Sciences, Fudan
University, Shanghai, 200433, P.R. China}\\
{\footnotesize \it E-mail: 051018003@fudan.edu.cn}}

\renewcommand{\thefootnote}{}
\maketitle

\pagestyle{myheadings} \thispagestyle{plain} \markboth{Lei Zhu
}{sharp lower bound}

\renewcommand{\baselinestretch}{1.5}
\begin{abstract} Let $V$ be a smooth projective 3-fold of general type.
Denote by $K^{3}$, a rational number, the self-intersection of the
canonical sheaf of any minimal model of $V$. One defines $K^{3}$ as
a canonical volume of $V$. The paper is devoted to proving the sharp
lower bound $K^{3}\ge \frac{1}{420}$ which can be reached by an
example: $X_{46}\subseteq \mathbb{P}(4,5,6,7,23)$.
\end{abstract}
\section{Introduction}
To classify algebraic varieties is one of the main goals of
algebraic geometry. Some explainations about the explicit structure
of algebraic threefolds can be found in the book \cite{10} edited by
Corti and Reid.

Let $V$ be a smooth projective 3-fold of general type. The
3-dimensional MMP allows us to pick a minimal model $X$ of $V$ and
$X$ has at worst $\mathbb{Q}$-factorial terminal singularities. A
recent work of Takayama (\cite{Tak}) and Hacon-M$^{\text{\rm
c}}$Kernan (\cite{H-M}) established the existence of a lower bound
for the canonical volume in all dimensions. There have been some
concrete known bounds already. In \cite{0}, J.A. Chen and M. Chen
proved $K_{X}^{3}\ge \frac{1}{30}$ for all 3-folds of general type
with $\chi(\Co{X})\le 0$. In \cite{11}, M. Chen gives a sharp lower
bound for the canonical volume of 3-folds of general type with
$P_{g}(V)\ge 2$.

In this paper, we study the sharp lower bound of volume of 3-folds
of general type with $\chi(\Co{X})=1$. The main idea is: firstly, we
study the pluricanonical maps of 3-folds and set a bound for
$P_{12}(X)$; then by analyzing the singularities of 3-folds of
general type and searching by some computable algorithm implemented
with Matlab, we obtained all the cases which may satisfy our
assumptions. Finally, we calculate $K_{X}^{3}$ and get the sharp
lower bound of the volumes.

The main result of this paper is the following:
\begin{maintheorem}
Let $X$ be a projective minimal 3-fold of general type with only
$\mathbb{Q}$-factorial terminal singularities. If $\chi(\Co{X})=1$,
then $K_{X}^{3}\ge \frac{1}{420}$.
\end{maintheorem}

\begin{example}
We found by our method one example which can reach the lower bound:
$\chi(\Co{X})=1$ and singularities: 3 of type $\frac{1}{2}(1,-1,1)$
and 1 of type $\frac{1}{4}(1,-1,1)$, $\frac{1}{5}(1,-1,2)$,
$\frac{1}{6}(1,-1,1)$ and $\frac{1}{7}(1,-1,3)$ and
$K_{X}^{3}=\frac{1}{420}$.  This is the hypersurface in weighted
projective space: $X_{46}\subseteq \mathbb{P}(4,5,6,7,23)$. One can
refer to the list of canonical 3-fold hypersurfaces (at page 358 of
\cite{12}) for details.
\end{example}

After completing this paper, the author was informed that the same
lower bound has also been obtained by J.A. Chen and M. Chen in
\cite[Theorem 6.11]{7}. We remark that these two approaches are
obviously different. This paper benefits less from so-called
theoretical skills, while more from the techniques of the
development of computable programming algorithm for
excluding most impossible cases.\\
{\bf Acknowledge.} This paper was prepared during my visit in Univ.
of Kent. I would like to thank Dr. G. Brown for invaluable helps and
hospitality. And I also thank Prof. Meng Chen for his guidance and
encouragement.

\section{Set up}
Suppose that $X$ is a projective minimal 3-fold of general type with
only $\mathbb{Q}$-factorial terminal singularities. If $P_{m}(X)\ge
2$, we take a resolution $\pi_{m}: X'\rightarrow X$, according to
Hironaka \cite{3}, such that:

\begin{enumerate}

\item[(i)] $X'$ is smooth;

\item[(ii)] the movable part of $|mK_{X'}|$ defines a morphism;

\item[(iii)] the fractional part of $\pi^*(mK_{X})$ has supports with only normal crossings. Denote by $g_{m}$ the
composition of $\Phi_{m}\circ\pi_{m}$. So $g_{m}:X'\rightarrow
W\subseteq \mathbb{P}^{P_{m}(X)-1}$ is a morphism. Let
$g_{m}:=s_{m}\circ f_{m}$ be the Stein factorization of $g_{m}$. We
can write $mK_{X'}\simlin M_{m}+Z_{m}$ and $mK_{X'}\simlin
\pi_{m}^*(mK_X)+E_m$ where $Z_{m}$ is the fixed part of $|mK_{X'}|$.
We also can write $\pi_{m}^*(mK_{X})\simlin M_{m}+\bar{E_{m}}$,
where $\bar{E_{m}}=Z_{m}-E_{m}$ is actually an effective
$\mathbb{Q}$-divisor. We have the following diagram:

\vspace{5mm}\hspace*{-1cm}\begin{picture}(50,80) \put(150,0){$X$}
\put(150,70){$X'$} \put(220,0){$W$} \put(220,70){$B$}
\put(162,75){\vector(1,0){53}} \put(154,65){\vector(0,-1){53}}
\put(225,65){\vector(0,-1){53}} \put(162,67){\vector(1,-1){53}}
\multiput(162,2.6)(5,0){11}{-} \put(212,5){\vector(1,0){5}}
\put(183,80){$f_{m}$} \put(230,35){$s_{m}$} \put(138,35){$\pi_{m}$}
\put(183,10){$\Phi_{m}$}\put(188,45){$g_{m}$}
\end{picture}
\end{enumerate}

If $\dim\Phi_{m}(X)=2$, we see that a general fiber of $f_{m}$ is a
smooth projective curve of genus $g\ge 2$. We say that $X$ is
canonically fibred by curves of genus $g$.

If $\dim\Phi_{m}(X)=1$, we see that a general fiber $S'_{m}$ of
$f_{m}$ is a smooth projective surface of general type. We say that
$X$ is canonically fibred by surfaces. We may write $M_{m}\simnum
aS'_{m}$ where $a\ge P_{m}(X)-1$. Particularly, if
$B$=$\mathbb{P}^{1}$, then we say that $|mK_{X}|$ is composed of a
rational pencil. Otherwise, we say that $|mK_{X}|$ is composed of a
irrational pencil.

Throughout this paper the symbol $\equiv$ stands for the nunmerical
equivalence of divisors, whereas $\sim_{lin}$ denotes the linear
equivalence and $=_{\mathbb{Q}}$ denotes the $\mathbb{Q}$-linear
equivalence.


\section{Some lemmas}
Following the techniques developed by M. Chen \cite{0,0000,chenzhu},
we can prove the next couple of lemmas which are used to show our
main theorem.

\begin{lemma}\label{2}
Let $X$ be a projective minimal 3-fold of general type with only
$\mathbb{Q}$-factorial terminal singularities. If $P_{12}(X)\ge 5$,
then $K_{X}^{3}\ge \frac{1}{360}$.
\end{lemma}
\begin{proof}
We analyze the bound in the following three cases.\\Case 1.
$\dim\phi_{12}(X)=3$.\\In this case, we can write
$12\pi_{12}^*(K_{X})|_{S_{12}}\simlin S_{12}+\bar{E_{12}} $ and
$|S_{12}|_{S_{12}}|$ is base point free and is not composed of a
pencil. Thus a generic irreducible element $C_{12}$ of
$|S_{12}|_{S_{12}}|$ is a smooth curve.
\begin{claim}The following holds
$$S_{12}|_{S_{12}}.C_{12}\ge 4.$$
\end{claim}
\begin{proof}
Let $D=(S_{12}|_{S_{12}})|_{C_{12}}$. Obviously, $\deg
D=S_{12}|_{S_{12}}.C_{12}$. Moreover, we have
\begin{eqnarray*}h^{0}(S_{12},S_{12}|_{S_{12}})\ge P_{12}(X)-1\ge 4,\\
h^{0}(C_{12},D)\ge h^{0}(S_{12},S_{12}|_{S_{12}})-1\ge
3,\end{eqnarray*} If $h^{1}(C_{12},D)=0$, noting that $C_{12}$ lives
in a covering family of a general type 3-fold and is hence of genus
at least $2$, then
\begin{eqnarray}h^{0}(C_{12},D)-h^{1}(C_{12},D)=h^{0}(C_{12},D)=1+\deg
D-g(C_{12})\end{eqnarray} which implies $\deg D\ge 4$.

If $h^{1}(C_{12},D)>0$, by Clifford theorem, one has
\begin{eqnarray}\deg D\ge 2h^{0}(C_{12},D)-2\ge 4.
\end{eqnarray}
Then, the desired result follows from  (1) and (2) at once.
\end{proof}
Therefore,
$2g(C_{12})-2=K_{S_{12}}.C_{12}+C_{12}^{2}=(K_{X'}+2S_{12})|_{S_{12}}.C_{12}\ge
10$, noting that $2g(C_{12})-2\equiv 0$(\text{mod} $2$). By theorem
2.4 in \cite{chenzhu}, we have
$m\pi_{12}^*(K_{X})|_{S_{12}}.C_{12}\ge
10+(m-25)\pi_{12}^*(K_{X})|_{S_{12}}.C_{12}$, which implies
$25\pi_{12}^*(K_{X})|_{S_{12}}.C_{12}\ge 10$ for $m>>0$. Therefore,
$12K_{X}^{3}\ge (\pi_{12}^*(K_{X}))^{2}.S_{12}\ge
\frac{1}{12}\pi_{12}^*(K_{X})|_{S_{12}}.C_{12}\ge \frac{10}{12\times
25}$, which means $K_{X}^{3}\ge \frac{1}{360}$.
\\Case 2.
$\dim\phi_{12}(X)=2$.
\\Recall that $S_{12}$ is a generic irreducible
element of $|M_{12}|$. $S_{12}|_{S_{12}}\simnum bC_{12}$ where $b\ge
h^{0}(S_{12},S_{12}|_{S_{12}})-1\ge P_{12}(X)-2\ge 3$.

If we take $m>>0$, then by theorem 2.4 in \cite{chenzhu}, we have
$m\pi_{12}^*(K_{X})|_{S_{12}}.C_{12}\ge
2g(C_{12})-2+(m-13-\frac{12}{b})\pi_{12}^*(K_{X})|_{S_{12}}.C_{12}$,
which implies $17\pi_{12}^*(K_{X})|_{S_{12}}.C_{12}\ge 2$. Next we
take $m_{1}=26$. We have $26\pi_{12}^*(K_{X})|_{S_{12}}.C_{12}\ge
2+\roundup{(26-13-\frac{12}{3})\times \frac{2}{17}}\ge 4$, which
means $\pi_{12}^*(K_{X})|_{S_{12}}.C_{12}\ge \frac{2}{13}$. Then
$$12K_{X}^{3}\ge \pi_{12}^*(K_{X})^{2}.S_{12}\ge
\frac{1}{12}\pi_{12}^*(K_{X})|_{S_{12}}.S_{12}|_{S_{12}}\ge
\frac{b}{12}\pi_{12}^*(K_{X})|_{S_{12}}.C_{12}\ge \frac{3}{12}\times
\frac{2}{13},$$ which implies $K_{X}^{3}\ge \frac{1}{312}$.\\Case 3.
$\dim\phi_{12}(X)=1$. \\In this case, we have
$12\pi_{12}^*(K_{X})\simnum aS'_{12}+E$ where $a\ge P_{12}(X)-1\ge
4$ and $S'_{12}$ is the general fiber. If $g(B)>0$, then $q(X)>0$.
J.A. Chen and Hacon have proved that $\phi_{8}$ is stably birational
(see \cite{C-H}). Therefore, we have $(8K_{X}^{3})\ge 2$ which
implies $K_{X}^{3}\ge \frac{1}{256}$.

If $g(B)=0$, according to Koll\'ar's technique (one can refer to
\cite{6}), we have $\mathcal{O}(4)\hookrightarrow
f_{m*}\omega_{X'}^{12}$ which means
$f_{m*}\omega_{X'/{\mathbb{P}^{1}}}^{2} \hookrightarrow
f_{m*}\omega_{X'}^{14}$. Therefore
$14\pi_{12}^*(K_{X})|_{S'_{12}}=_{\mathbb{Q}}
2\sigma^*(K_{S_{0}})+F_{12}$ where $F_{12}$ is an effective
${\mathbb{Q}}$-divisor. This implies
$$(7\pi_{12}^*(K_{X})|_{S'_{12}})^{2}=49(\pi_{12}^*(K_{X})|_{S'_{12}})^{2}\ge
1.$$ Therefore, $12K_{X}^{3}\ge \pi_{12}^*(K_{X})^{2}.(aS'_{12})\ge
4(\pi_{12}^*(K_{X})|_{S'_{12}})^{2}\ge \frac{4}{49}$.
\end{proof}

\begin{lemma}\label{3}
Let $X$ be a projective minimal 3-fold of general type with only
$\mathbb{Q}$-factorial terminal singularities. If $P_{m}(X)=2$, then
$K_{X}^{3}\ge \frac{1}{m(m+1)^{2}}$.
\end{lemma}
\begin{proof}
Since $P_{m}(X)=2$, $\dim\phi_{m}(X)=1$. And we can write
$m\pi_{m}^*(K_{X})\simnum aS'_{m}+\bar{E_{m}}$ where $a\ge
P_{m}(X)-1=1$.\\If $g(B)>0$, then $K_{X}^{3}\ge \frac{1}{256}$ for
the same reason as above.\\If $g(B)=0$, due to Koll\'ar's technique,
we have $\mathcal{O}(1)\hookrightarrow f_{m*}\omega_{X'}^{m}$ which
implies
\\$\pi_{m}^*(K_{X})|_{S'_{m}}\simnum
\beta_{n}\sigma^*(K_{S_{0}})+F_{m}$ due to lemma 2.8 in
\cite{chenzhu}, where $\lim_{n\mapsto +\infty} \beta_n =
\frac{1}{m+1}$ and $F_{m}$ is an effective $\mathbb{Q}$-divisor.

If $P_{g}(S'_{m})>0$, one can take $C_{2}$ the general member of
$|2\sigma^*(K_{S_{0}})|$ since $|2\sigma^*(K_{S_{0}})|$ is base
point free by \cite{8}. Thus
$(2m+2)\pi_{m}^*(K_{X})|_{S'_{m}}.C_{2}\ge 4$. Then $mK_{X}^{3}\ge
\frac{1}{m+1}\pi_{m}^*(K_{X})|_{S'_{m}}.C_{2}\ge
\frac{1}{(m+1)^{2}}$ which says $K_{X}^{3}\ge
\frac{1}{m(m+1)^{2}}$.

If $P_{g}(S'_{m})=0$, $|3\sigma^*(K_{S_{0}})|$ has no fixed part by
[Bombieri, Reider, Miyaoka, Catanese]. Take $C_{3}$ the general
member of $|3\sigma^*(K_{S_{0}})|$. Thus
$$(3m+3)\pi_{m}^*(K_{X})|_{S'_{m}}.C_{3}\ge 9.$$ Therefore,
$mK_{X}^{3}\ge \frac{1}{(3m+3)}\pi^*(K_{X})|_{S'_{m}}.C_{3}\ge
\frac{1}{(m+1)^{2}}$ which gives $K_{X}^{3}\ge
\frac{1}{m(m+1)^{2}}$.
\end{proof}

\section{The proof of Main Theorem 1}
We always assume $X$ to be a projective minimal 3-fold of general
type with only $\mathbb{Q}$-factorial terminal singularities and
with $\chi(\Co{X})=1$. Let us recall Reid's plurigenera formula (at
page 413 of \cite{12}) for a minimal 3-fold $X$ of general type:
$$P_{m}(X)=(1-2m)\chi(\mathcal{O}_X)+\frac{1}{12}m(m-1)(2m-1)K_{X}^{3}+l(m),$$
where $m>1$ is an integer, the correction term
$$l(m):=\sum_{Q}l(Q,m):=\sum_{Q}\sum_{j=1}^{m-1}\frac{\bar{bj}(r-\bar{bj})}{2r},$$
where the sum $\sum_{Q}$ runs through all baskets $Q$ of
singularities of type $\frac{1}{r}(a,-a,1)$ with the positive
integer $a$ coprime to $r$, $0<a<r$, $0<b<r$, $ab\equiv 1$
(\text{mod} $r$), $\bar{bj}$ the smallest residue of $bj$ mod $r$.
Reid's result (Theorem 10.2 in \cite{12}) says that the above
baskets ${Q}$ of singularities are in fact virtual (!) and that one
need not worry about the authentic type of all those terminal
singularities on $X$, though $X$ may have nonquotient terminal
singularities. Iano-Fletcher (\cite{Fletcher}) has showed that the
set of baskets ${Q}$ in Reid's formula is uniquely determined by
$X$.

\begin{lemma}\cite[Lemma 3.1]{2}\label{4}
For all $m\ge 0$ and $n\ge 1$, we have$$l(m+2n)\ge l(m)+nl(2)$$ with
equality if and only if all the singularities of type
$\frac{1}{2}(1,1,1)$.
\end{lemma}

\begin{lemma}\cite[lemma 3.2]{2}\label{5}
Suppose $\alpha>\beta$ are integers. Then we have
$$l(\frac{1}{\alpha}(1,-1,1),n)\ge l(\frac{1}{\beta}(1,-1,1),n)$$
for all $n<\beta<\alpha$.
\end{lemma}

\begin{lemma}\cite[lemma 3.3]{2}\label{6}
For all $n\le [\frac{r+1}{2}]$, we have
$$l(\frac{1}{r}(a,-a,1),n)\ge l(\frac{1}{r}(1,-1,1),n),$$
where $[\frac{r+1}{2}]$ denotes the integral part of
$\frac{r+1}{2}$.
\end{lemma}

\begin{lemma}\cite[corollary 3.4]{2}\label{7}
For all $\alpha$, $\beta\in \mathbf{Z}$ with $0\le \beta \le \alpha$
and  $n\le [\frac{\alpha+1}{2}]$, we have
$$l(\frac{1}{\alpha}(a,-a,1),n)\ge l(\frac{1}{\beta}(1,-1,1),n).$$
\end{lemma}

\begin{lemma}\label{8}
If $P_{2}(X)\ge 1$, then $P_{2n}(X)\ge n$.
\end{lemma}
\begin{proof}
By Lemma \ref{4}, we have
$$P_{4}(X)=-7+7K_{X}^{3}+l(4) \quad \mbox{ and }\quad l(4)\ge 2l(2).$$
Since $$P_{2}(X)=-3+\frac{1}{2}K_{X}^{3}+l(2),$$we get
$$P_{4}(X)-2P_{2}(X)=-1+6K_{X}^{3}+l(4)-2l(2).$$ This implies $$P_{4}(X)\ge
1+6K_{X}^{3}.$$ Thus we have
\begin{eqnarray*}P_{4}(X)&\ge& 2,\\
P_{6}(X)&=&-11+\frac{55}{2}K_{X}^{3}+l(6),\\
P_{6}(X)-P_{4}(X)-P_{2}(X)&=&-11+7+3+(\frac{55}{2}-7-\frac{1}{2})K_{X}^{3}+l(6)-l(4)-l(2)\\&\ge&
-1+20K_{X}^{3}.
\end{eqnarray*}

It follows from Lemma \ref{4} that  $l(6)\ge l(4)+l(2)$. Simple
calculation yields
$$P_{6}(X)\ge 2+20K_{X}^{3}>2.$$

Assume $P_{2m}(X)\ge m$, then we have
$$P_{2m+2}(X)=(-3-4m)+\frac{(2m+2)(2m+1)(4m+3)}{12}K_{X}^{3}+l(2m+2)$$
and
$$P_{2m+2}(X)-P_{2}(X)-P_{2m}(X)=-1+kK_{X}^{3}+l(2m+2)-l(2m)-l(2),$$
where $k>0$. Then $P_{2m+2}(X)\ge m+kK_{X}^{3}$ which means
$P_{2m+2}(X)\ge m+1$.
\end{proof}

\begin{lemma}\label{9}
If $P_{6}(X)\ge 3$, then $P_{12}(X)\ge 5$.
\end{lemma}
\begin{proof}If $P_{6}(X)\ge 4$, we have $P_{12}(X)=h^{0}(X',12K_{X'})>P_{6}(X)\ge 5$.
Otherwise, $6K_{X'}$ is fixed which is a contradiction. Assume
$P_{6}(X)=3$ and $P_{12}(X)=4$. Keep the same notation as in section
2. We study $f_{6}: X'\rightarrow \mathbb{P}^{2}$. If $|6K_{X'}|$ is
composed of pencils, then $P_{12}(X)\ge 5$, a contradiction.
Otherwise, we can write $|12K_{X'}|\supset
|K_{X'}+\roundup{5\pi_{6}^*(K_{X})}+S_{6}|$, where $S_{6}$ is the
general member of the movable part of $|6K_{X'}|$. Consider the
following exact sequence: $$0\rightarrow
H^0(X',K_{X'}+\roundup{5\pi_{6}^*(K_{X})})\rightarrow
H^0(X',K_{X'}+\roundup{5\pi_{6}^*(K_{X})}+S_{6})$$ $$\rightarrow
H^0(S_{6},K_{S_{6}}+\roundup{5\pi_{6}^*(K_{X})}|_{S_{6}})\rightarrow
0,$$ the surjective map is due to Kawamata-Viehweg vanishing
theorem. Therefore,
$$
P_{12}(X)\ge
h^0(X',K_{X'}+\roundup{5\pi_{6}^*(K_{X})}+S_{6})=P_{6}+h^0(S_{6},K_{S_{6}}+\roundup{5\pi_{6}^*(K_{X})}|_{S_{6}})\ge
5,$$ noting that
$$h^0(S_{6},K_{S_{6}}+\roundup{5\pi_{6}^*(K_{X})}|_{S_{6}})=h^0(S_{6},(K_{X'}+S_{6}+\roundup{5\pi_{6}^*(K_{X})})|_{S_{6}})\ge
2.$$ So we are done.
\end{proof}

\begin{proof}[{\bf Proof of Main Theorem.}]
If $P_{6}(X)\ge 2$ or $P_{4}(X)\ge 2$, one can conclude that
$K_{X}^{3}>\frac{1}{420}$ by lemma \ref{2}, lemma \ref{3} and lemma
\ref{9}.

Firstly, we deal with the case $P_{12}(X)=4$. Due to Lemma \ref{8},
we have $P_{2}(X)=0$. If $P_{3}(X)\ge 2$, then $P_{6}(X)\ge 3$.
Otherwise, one can derive $3K_{X'}$ is fixed which is a
contradiction. Thus $P_{12}(X)\ge 5$. Since $P_{6}(X)\le 1$ and
$P_{4}(X)\le 1$, we have the following five cases:
\begin{enumerate}

\item[(i)] $P_{2}(X)=0$, $P_{3}(X)=0$, $P_{4}(X)=0$ and
$P_{6}(X)=0$;

\item[(ii)] $P_{2}(X)=0$, $P_{3}(X)=0$, $P_{4}(X)=0$ and
$P_{6}(X)=1$;

\item[(iii)] $P_{2}(X)=0$, $P_{3}(X)=0$, $P_{4}(X)=1$ and
$P_{6}(X)=1$;

\item[(iv)] $P_{2}(X)=0$, $P_{3}(X)=1$, $P_{4}(X)=0$ and
$P_{6}(X)=1$;

\item[(v)] $P_{2}(X)=0$, $P_{3}(X)=1$, $P_{4}(X)=1$ and
$P_{6}(X)=1$.
\end{enumerate}

We use the matrix $E$ (firstly used by Fletcher \cite{2} and such
matrix is not unique)
$$ \left(
\begin{array}{cccc}
3 &-2 &-4 &-10\\
-1 &1 &-3 &-16\\
0 &0 &1 &-7\\
0 &0 &0 &1\\
\end{array}
\right)
$$
Obviously, one has the following equalities due to Reid's formula
(at page 413 in \cite{12}), setting
$\Delta_{n}:=n^{2}l(2)+l(n)-l(n+1)$:
$$
\left\{
\begin{array}{l}
P_{3}(X)=10-\sum_{Q}\Delta_{2}(Q)\\
P_{4}(X)-P_{3}(X)=25-\sum_{Q}\Delta_{3}(Q)\\
P_{6}(X)-P_{4}(X)=119-\sum_{Q}(\Delta_{4}(Q)+\Delta_{5}(Q))\\
P_{12}(X)-P_{6}(X)=1341-\sum_{Q}(\Delta_{6}(Q)+...+\Delta_{11}(Q))
\end{array}
\right.
$$
For brief, we denote by $\nabla_{1}=\Delta_{2}$,
$\nabla_{2}=\Delta_{3}$, $\nabla_{3}=\Delta_{4}+\Delta_{5}$ and
$\nabla_{4}=\Delta_{6}+...+\Delta_{11}$.

\begin{claim}
If $X$ contains a singularity of type $\frac{1}{s}(a, -a, 1)$ where
$s\ge 43$, then $P_{12}(X)\ge 5$.
\end{claim}
\begin{proof}
Let $Q$ be a singularity of type $\frac{1}{s}(a, -a, 1)$ and $Q'$ be
a singularity of type $\frac{1}{s}(1, -1, 1)$. Then by lemma
\ref{7}, $$l(12)\ge l(Q, 12)\ge l(Q',
12)=\sum_{j=1}^{11}\frac{j\times (43-j)}{2\times 43}\ge 27.$$ On the
other hand, due to Reid \cite{12},
$$P_{12}(X)=-23+253K_{X}^{3}+l(12)> -23+27,$$ which means $P_{12}(X)\ge
5$.
\end{proof}

All the datum are listed in Appendix 1.

For Case (i). We have
$$
\left\{
\begin{array}{l}
\sum_{Q}\nabla_{1}(Q)=10\\
\sum_{Q}\nabla_{2}(Q)=25\\
\sum_{Q}\nabla_{3}(Q)=119\\
\sum_{Q}\nabla_{4}(Q)=1337
\end{array}
\right.
$$
Reducing $\nabla$ via the matrix $E$, gives four new equations:
$$
\left\{
\begin{array}{l}
\sum_{Q}\nabla'_{1}(Q)=5\\
\sum_{Q}\nabla'_{2}(Q)=5\\
\sum_{Q}\nabla'_{3}(Q)=4\\
\sum_{Q}\nabla'_{4}(Q)=4
\end{array}
\right.
$$
By reference to Appendix 1, one solution is found: $2\times
\frac{1}{5}(3,-3,1)$, $3\times \frac{1}{2}(1,-1,1)$, $2\times
\frac{1}{3}(2,-2,1)$ and $\frac{1}{4}(3,-3,1)$. Now
$l(2)=\frac{359}{120}$ which implies $K_{X}^{3}=\frac{1}{60}$.

In the case (ii) through (v), one has, respectively:
\begin{eqnarray*}
&&(\nabla_1',\nabla_2',\nabla_3',\nabla_4')=
(\nabla_1,\nabla_2,\nabla_3,\nabla_4)E\\
&=&(5,5,3,12), (6,4,7,21), (1,8,4,6),(2,7,8,15).
\end{eqnarray*}
Using the same technique as above and searched by computer, the all
solutions of case 2 to case 5 respectively are:
\renewcommand{\arraystretch}{1.5}
$$
\begin{array}{cclc}\hline
Case. \hspace*{2mm}&\hspace*{2mm} P_{12}(X) \hspace*{2mm}
&\hspace*{20mm} Singularities \hspace*{2mm}
&\hspace*{2mm}K_{X}^{3}\\ ii. &4  &\frac{1}{7}(4,-4,1),
\frac{1}{9}(7,-7,1), 2\times
\frac{1}{2}(1,-1,1), 2\times \frac{1}{3}(2,-2,1)& \frac{1}{63}\\
iii. &4& \frac{1}{9}(5,-5,1), \frac{1}{16}(9,-9,1),
\frac{1}{2}(1,-1,1)&\frac{1}{144}\\
iii. &4 &2\times \frac{1}{5}(4,-4,1), \frac{1}{9}(7,-7,1), 3\times
\frac{1}{2}(1,-1,1), \frac{1}{3}(2,-2,1)&\frac{1}{90}\\
iii.& 4 & \frac{1}{7}(5,-5,1),\frac{1}{9}(7,-7,1),
\frac{1}{9}(5,-5,1), \frac{1}{2}(1,-1,1)&\frac{1}{126}\\
iii.&4& \frac{1}{4}(3,-3,1), \frac{1}{5}(4,-4,1),
\frac{1}{7}(5,-5,1), \frac{1}{11}(9,-9,1)&\frac{13}{1540}\\
v. &4 &3\times \frac{1}{4}(3,-3,1), \frac{1}{5}(3,-3,1),
\frac{1}{11}(8,-8,1)&\frac{1}{220}\\
v. &4&2\times \frac{1}{4}(3,-3,1), 2\times \frac{1}{5}(3,-3,1),
\frac{1}{7}(4,-4,1), \frac{1}{3}(2,-2,1)&\frac{1}{210}\\ \hline
\end{array}
$$

As for $P_{12}=3$, we have the same five cases as above and the corresponding datum:
\begin{eqnarray*}
(\nabla'_1,\nabla'_2,\nabla'_3,\nabla'_4)& =&(5,5,4,5),
(5,5,3,13),\\
&& (6,4,7,22), (1,8,4,7), (2,7,8,16).
\end{eqnarray*}

Then one can get the following table searched by
computer:

$$
\begin{array}{cclc}\hline
Case. \hspace*{2mm}&\hspace*{2mm} P_{12}(X) \hspace*{2mm} &\hspace*{20mm} Singularities \hspace*{2mm} &\hspace*{2mm}K_{X}^{3}\\
i. &3& \frac{1}{7}(5,-5,1), \frac{1}{18}(13,-13,1) & \frac{1}{126}\\
i. &3&\frac{1}{8}(5,-5,1), \frac{1}{17}(12,-12,1) & \frac{1}{136}\\
i. &3&\frac{1}{12}(7,-7,1), \frac{1}{13}(8,-8,1) &\frac{1}{156}\\
i. &3&\frac{1}{5}(3,-3,1), \frac{1}{7}(5,-5,1),
\frac{1}{13}(8,-8,1)&\frac{4}{455}\\
i. &3&\frac{1}{5}(3,-3,1), \frac{1}{8}(5,-5,1), \frac{1}{12}(7,-7,1)&
\frac{1}{120}\\
i. &3&2\times \frac{1}{5}(3,-3,1), \frac{1}{7}(5,-5,1),
\frac{1}{8}(5,-5,1)&\frac{3}{280}\\
ii. &3&\frac{1}{7}(4,-4,1), \frac{1}{11}(9,-9,1),
\frac{1}{2}(1,-1,1),
2\times \frac{1}{3}(2,-2,1)&\frac{5}{462}\\
ii. &3&\frac{1}{7}(4,-4,1), \frac{1}{13}(11,-11,1), 2\times
\frac{1}{3}(2,-2,1)&\frac{2}{273}\\
ii. &3& \frac{1}{9}(7,-7,1), \frac{1}{10}(7,-7,1),
\frac{1}{3}(2,-2,1), 2\times \frac{1}{2}(1,-1,1)&\frac{1}{90}\\
ii. &3& \frac{1}{9}(7,-7,1), \frac{1}{13}(10,-10,1), 2\times
\frac{1}{2}(1,-1,1)&\frac{1}{117}\\
iii. &3&2\times \frac{1}{5}(4,-4,1), \frac{1}{11}(9,-9,1), 2\times
\frac{1}{2}(1,-1,1), \frac{1}{3}(2,-2,1)&\frac{1}{165}\\
iii. &3 & 2\times \frac{1}{5}(4,-4,1), \frac{1}{13}(11,-11,1),
\frac{1}{2}(1,-1,1), \frac{1}{3}(2,-2,1)&\frac{1}{390}\\
iii.&3&\frac{1}{7}(5,-5,1), \frac{1}{9}(5,-5,1),
\frac{1}{11}(9,-9,1)&\frac{2}{693}\\
iii. &3&\frac{1}{4}(3,-3,1), \frac{1}{5}(4,-4,1), 2\times
\frac{1}{9}(7,-7,1)&\frac{1}{180}\\
iii. &3&\frac{1}{4}(3,-3,1), \frac{1}{5}(3,-3,1),
\frac{1}{6}(5,-5,1),
\frac{1}{7}(5,-5,1), 3\times \frac{1}{2}(1,-1,1)&\frac{1}{420}\\
v. &3&\frac{1}{4}(3,-3,1),2\times \frac{1}{5}(3,-3,1),
\frac{1}{11}(7,-7,1), \frac{1}{3}(2,-2,1)&\frac{1}{660}\\\hline
\end{array}
$$

The last but two solution is the codimension 1 hypersurface
$X_{46}\subset \mathbb{P}(4,5,6,7,23)$ in \cite{10}.

As for $P_{12}(X)=1$ and $P_{12}(X)=2$, we have the same cases.
Noting that $P_{4}(X)\le 1$ and $P_{6}(X)\le 1$ in these two cases, our
table is integrated:
$$
\begin{array}{cclc}\hline
Case. \hspace*{2mm}&\hspace*{2mm} P_{12}(X) \hspace*{2mm} &\hspace*{20mm} Singularities \hspace*{2mm} &\hspace*{2mm}K_{X}^{3}\\
i. &1&\frac{1}{4}(3,-3,1), \frac{1}{5}(3,-3,1), \frac{1}{7}(5,-5,1),
2\times \frac{1}{2}(1,-1,1), 2\times
\frac{1}{3}(2,-2,1)&\frac{1}{420}\\
i. &1&\frac{1}{4}(3,-3,1), \frac{1}{5}(3,-3,1), \frac{1}{11}(8,-8,1),
3\times\frac{1}{2}(1,-1,1)&\frac{1}{220}\\
i. &1&2\times \frac{1}{5}(3,-3,1), \frac{1}{7}(4,-4,1), 3\times
\frac{1}{2}(1,-1,1), \frac{1}{3}(2,-2,1)&\frac{1}{210}\\
i. &2&\frac{1}{4}(3,-3,1), \frac{1}{13}(8,-8,1), 3\times
\frac{1}{2}(1,-1,1), \frac{1}{3}(2,-2,1)&\frac{1}{156}\\
i. &2&\frac{1}{4}(3,-3,1), \frac{1}{5}(3,-3,1), \frac{1}{8}(5,-5,1),
\frac{1}{3}(2,-2,1), 3\times \frac{1}{2}(1,-1,1)&\frac{1}{120}\\
ii. &2&\frac{1}{10}(7,-7,1), \frac{1}{11}(9,-9,1),
\frac{1}{2}(1,-1,1), \frac{1}{3}(2,-2,1)&\frac{1}{165}\\
ii.&2&\frac{1}{10}(7,-7,1), \frac{1}{13}(11,-11,1),
\frac{1}{3}(2,-2,1)&\frac{1}{390}\\
ii. &2&\frac{1}{11}(9,-9,1), \frac{1}{13}(10,-10,1),
\frac{1}{2}(1,-1,1)&\frac{1}{286}\\\hline
\end{array}
$$

\begin{claim}
The solution $\frac{1}{4}(3,-3,1)$,2 of $\frac{1}{5}(3,-3,1)$,
$\frac{1}{11}(7,-7,1)$ and $\frac{1}{3}(2,-2,1)$ and
$K_{X}^{3}=\frac{1}{660}$ does not exist.
\end{claim}
\begin{proof}
Since $P_{7}(X)=-13+\frac{7\times 6\times 13}{12}K_{X}^{3}+l(7)=2$,
one can derive $K_{X}^{3}\ge \frac{1}{448}$ by lemma \ref{3}, a
contradiction.
\end{proof}

\begin{claim}
The solution $\frac{1}{4}(3,-3,1)$, $\frac{1}{5}(3,-3,1)$,
$\frac{1}{7}(5,-5,1)$, $2\times \frac{1}{2}(1,-1,1)$, $2\times
\frac{1}{3}(2,-2,1)$ and $K_{X}^{3}=\frac{1}{420}$ does not exist.
\end{claim}
\begin{proof}
According to Reid (see (10.3) of \cite{12}), one has
$$\frac{1}{12}K_X\cdot c_2=-2\chi({\mathcal
O}_X)+\sum_Q\frac{r_Q^2-1}{12r_Q}.$$ Miyaoka \cite{133} says
$K_X\cdot c_2\geq 0$. Thus we get the inequality
$$\sum_Q\frac{r_Q^2-1}{r_Q}\geq 24\chi({\mathcal
O}_X)=24.$$ Now since the datum in this claim doesn't fit into the
above inequality, this solution doesn't exist. By the way, the
effectivity of Miyaoka-Reid inequality was firstly observed in
\cite{shin}.
\end{proof}
\end{proof}
\clearpage

\noindent {\bf Appendix 1.}
\renewcommand{\arraystretch}{1.5}
The following table gives the values of $\nabla_{i}$ and
$\nabla'_{i}$.
$$
\begin{array}{cccccccccc}\hline
No. \hspace*{2mm}&\hspace*{2mm}Singularity\hspace*{2mm} &\hspace*{2mm}\nabla_{1}\hspace*{2mm} &\hspace*{2mm}\nabla_{2}\hspace*{2mm} &\hspace*{2mm}\nabla_{3}\hspace*{2mm} \ &\hspace*{2mm}\nabla_{4}\hspace*{2mm} &\hspace*{2mm}\nabla'_{1} \hspace*{2mm}&\hspace*{2mm}\nabla'_{2} \hspace*{2mm}&\nabla'_{3} \hspace*{2mm}&\hspace*{2mm}\nabla'_{4}\\
1.  &\frac{1}{2}(1, -1, 1)    &1  &  2&     10\ &  112 & 1 & 0 & 0& 0\\
2.  &\frac{1}{3}(2, -2, 1)    &1  &  3&     13\ &  149 & 0 & 1 & 0& 0\\
3.  &\frac{1}{4}(3, -3, 1)    &1  &  3&     15\ &  167 & 0 & 1 & 2& 4\\
4.  &\frac{1}{5}(4, -4, 1)    &1  &  3&     16\ &  178 & 0 & 1 & 3& 8\\
5.  &\frac{1}{5}(3, -3, 1)    &2  &  5&     24\ &  268 & 1 & 1 & 1& 0\\
6.  &\frac{1}{6}(5, -5, 1)    &1  &  3&     16\ &  185 & 0 & 1 & 3& 15\\
7.  &\frac{1}{7}(6, -6, 1)    &1  &  3&     16\ &  190 & 0 & 1 & 3& 20\\
8.  &\frac{1}{7}(5, -5, 1)    &3  &  7&     34\ &  383 & 2 & 1 & 1& 3\\
9.  &\frac{1}{7}(4, -4, 1)    &2  &  6&     28\ &  319 & 0 & 2 & 2& 7\\
10.  &\frac{1}{8}(7, -7, 1)    &1  &  3&     16\ &  194 & 0 & 1 & 3& 24\\
11.  &\frac{1}{8}(5, -5, 1)    &3  &  8&     37\ &  419 & 1 & 2 & 1& 2\\
12.  &\frac{1}{9}(8, -8, 1)    &1  &  3&     16\ &  197 & 0 & 1 & 3& 27\\
13.  &\frac{1}{9}(7, -7, 1)    &4  &  9&     44\ &  497 & 3 & 1 & 1& 5\\
14.  &\frac{1}{9}(5, -5, 1)    &2  &  6&     31\ &  346 & 0 & 2 & 5& 13\\
15.  &\frac{1}{10}(9, -9, 1)    &1  &  3&     16\ &  199 & 0 & 1 & 3& 29\\
16.  &\frac{1}{10}(7, -7, 1)    &3  &  9&     41\ &  469 & 0 & 3 & 2& 8\\
17.  &\frac{1}{11}(10, -10, 1)    &1  &  3&     16\ &  200 & 0 & 1 & 3& 30\\
18.  &\frac{1}{11}(9, -9, 1)    &5  &  11&     54\ &  610 & 4 & 1 & 1& 6\\
19.  &\frac{1}{11}(8, -8, 1)    &4  &  11&     50\ &  569 & 1 & 3 & 1& 3\\
20.  &\frac{1}{11}(7, -7, 1)    &3  &  9&     43\ &  487 & 0 & 3 & 4& 12\\
21.  &\frac{1}{11}(6, -6, 1)    &2  &  6&     32\ &  364 & 0 & 2 &
6& 24\\\hline
\end{array}
$$
$$
\begin{array}{cccccccccc}\hline
No. \hspace*{2mm}&\hspace*{2mm}Singularity\hspace*{2mm} &\hspace*{2mm}\nabla_{1}\hspace*{2mm} &\hspace*{2mm}\nabla_{2}\hspace*{2mm} &\hspace*{2mm}\nabla_{3}\hspace*{2mm} \ &\hspace*{2mm}\nabla_{4}\hspace*{2mm} &\hspace*{2mm}\nabla'_{1} \hspace*{2mm}&\hspace*{2mm}\nabla'_{2} \hspace*{2mm}&\nabla'_{3} \hspace*{2mm}&\hspace*{2mm}\nabla'_{4}\\
22.  &\frac{1}{12}(11, -11, 1)    &1  &  3&     16\ &  200 & 0 & 1 & 3& 30\\
23.  &\frac{1}{12}(7, -7, 1)    &5  &  12&     58\ &  651 & 3 & 2 & 2& 3\\
24.  &\frac{1}{13}(12, -12, 1)    &1  &  3&     16\ &  200 & 0 & 1 & 3& 30\\
25.  &\frac{1}{13}(11, -11, 1)    &6  &  13&     64\ &  722 & 5 & 1 & 1& 6\\
26.  &\frac{1}{13}(10, -10, 1)    &4  &  12&     54\ &  618 & 0 & 4 & 2& 8\\
27.  &\frac{1}{13}(9, -9, 1)    &3  &  9&     46\ &  513 & 0 & 3 & 7& 17\\
28.  &\frac{1}{13}(8, -8, 1)    &5  &  13&     61\ &  687 & 2 & 3 & 2& 2\\
29.  &\frac{1}{13}(7, -7, 1)    &2  &  6&     32\ &  375 & 0 & 2 & 6& 35\\
30.  &\frac{1}{14}(13, -13, 1)    &1  &  3&     16\ &  200 & 0 & 1 & 3& 30\\
31.  &\frac{1}{14}(11, -11, 1)    &5 &  14&     63\ &  718 & 1 & 4 & 1& 3\\
32.  &\frac{1}{14}(9, -9, 1)    &3  &  9&     47\ &  524 & 0 & 3 & 8& 21\\
33.  &\frac{1}{15}(14, -14, 1)    &1  &  3&     16\ &  200 & 0 & 1 & 3& 30\\
34.  &\frac{1}{15}(13, -13, 1)    &7  &  15&     74\ &  834 & 6 & 1 & 1& 6\\
35.  &\frac{1}{15}(11, -11, 1)    &4  &  12&     58\ &  654 & 0 & 4 & 6& 16\\
36.  &\frac{1}{15}(8, -8, 1)    &2  &  6&     32\ &  384 & 0 & 2 & 6& 44\\
37.  &\frac{1}{16}(15, -15, 1)    &1  &  3&     16\ &  200 & 0 & 1 & 3& 30\\
38.  &\frac{1}{16}(13, -13, 1)    &5  &  15&     67\ &  767 & 0 & 5 & 2& 8\\
39.  &\frac{1}{16}(11, -11, 1)    &3  &  9&     48\ &  542 & 0 & 3 & 9& 32\\
40.  &\frac{1}{16}(9, -9, 1)    &7  &  16&     78\ &  880 & 5 & 2 & 2& 8\\
41.  &\frac{1}{17}(16, -16, 1)    &1  &  3&     16\ &  200 & 0 & 1 & 3& 30\\
42.  &\frac{1}{17}(15, -15, 1)    &8  &  17&    84\ &  946 & 7 & 1 &
1& 6\\\hline
\end{array}
$$
$$
\begin{array}{cccccccccc}\hline
No. \hspace*{2mm}&\hspace*{2mm}Singularity\hspace*{2mm} &\hspace*{2mm}\nabla_{1}\hspace*{2mm} &\hspace*{2mm}\nabla_{2}\hspace*{2mm} &\hspace*{2mm}\nabla_{3}\hspace*{2mm} \ &\hspace*{2mm}\nabla_{4}\hspace*{2mm} &\hspace*{2mm}\nabla'_{1} \hspace*{2mm}&\hspace*{2mm}\nabla'_{2} \hspace*{2mm}&\nabla'_{3} \hspace*{2mm}&\hspace*{2mm}\nabla'_{4}\\
43.  &\frac{1}{17}(14, -14, 1)    &6  &  17&     76\ &  867 & 1 & 5 & 1& 3\\
44.  &\frac{1}{17}(13, -13, 1)    &4  &  12&     61\ &  680 & 0 & 4 & 9& 21\\
45.  &\frac{1}{17}(12, -12, 1)    &7  &  17&     82\ &  919 & 4 & 3 & 3& 3\\
46.  &\frac{1}{17}(11, -11, 1)    &3  &  9&     48\ &  549 & 0 & 3 & 9& 39\\
47.  &\frac{1}{17}(10, -10, 1)    &5  &  15&     69\ &  788 & 0 & 5 & 4& 15\\
48.  &\frac{1}{17}(9, -9, 1)    &2  &  6&     32\ &  391 & 0 & 2 & 6& 51\\
49.  &\frac{1}{18}(17, -17, 1)    &1  &  3&     16\ &  200 & 0 & 1 & 3& 30\\
50.  &\frac{1}{18}(13, -13, 1)    &7  &  8&     85\ &  955 & 3 & 4 & 3& 2\\
51.  &\frac{1}{18}(11, -11, 1)    &5  &  15&     71\ &  806 & 0 & 5 & 6& 19\\
52.  &\frac{1}{19}(18, -18, 1)    &1  &  3&     16\ &  200 & 0 & 1 & 3& 30\\
53.  &\frac{1}{19}(17, -17, 1)    &9  &  19&     94\ &  1058 & 8 & 1 & 1& 6\\
54.  &\frac{1}{19}(16, -16, 1)    &6  &  18&     80\ &  916 & 0 & 6 & 2& 8\\
55.  &\frac{1}{19}(15, -15, 1)    &5  &  15&     73\ &  821 & 0 & 5 & 8& 20\\
56.  &\frac{1}{19}(14, -14, 1)    &4  &  12&     63\ &  702 & 0 & 4 & 11& 29\\
57.  &\frac{1}{19}(13, -13, 1)    &3  &  9&     48\ &  560 & 0 & 3 & 9& 50\\
58.  &\frac{1}{19}(12, -12, 1)    &8  &  19&     92\ &  1034 & 5 & 3 & 3& 6\\
59.  &\frac{1}{19}(11, -11, 1)    &7  &  19&     87\ &  988 & 2 & 5 & 2& 5\\
60.  &\frac{1}{19}(10, -10, 1)    &2  &  6&     32\ &  396 & 0 & 2 & 6& 56\\
61.  &\frac{1}{20}(19, -19, 1)    &1  &  3&     16\ &  200 & 0 & 1 & 3& 30\\
62.  &\frac{1}{20}(17, -17, 1)    &7  &  20&     89\ &  1016 & 1 & 6 & 1& 3\\
63.  &\frac{1}{20}(13, -13, 1)    &3  &  9&     48\ &  565 & 0 & 3 &
9& 55\\\hline
\end{array}
$$

$$
\begin{array}{cccccccccc}\hline
No. \hspace*{2mm}&\hspace*{2mm}Singularity\hspace*{2mm} &\hspace*{2mm}\nabla_{1}\hspace*{2mm} &\hspace*{2mm}\nabla_{2}\hspace*{2mm} &\hspace*{2mm}\nabla_{3}\hspace*{2mm} \ &\hspace*{2mm}\nabla_{4}\hspace*{2mm} &\hspace*{2mm}\nabla'_{1} \hspace*{2mm}&\hspace*{2mm}\nabla'_{2} \hspace*{2mm}&\nabla'_{3} \hspace*{2mm}&\hspace*{2mm}\nabla'_{4}\\
64.  &\frac{1}{20}(11, -11, 1)    &9  &  20&     98\ &  1107 & 7 & 2 & 2& 11\\
65.  &\frac{1}{21}(20, -20, 1)    &1  &  3&     16\ &  200 & 0 & 1 & 3& 30\\
66.  &\frac{1}{21}(19, -19, 1)    &10  &  21&     104\ &  1170 & 9 & 1 & 1& 6\\
67.  &\frac{1}{21}(17, -17, 1)    &5  &  15&     76\ &  847 & 0 & 5 & 11& 25\\
68.  &\frac{1}{21}(16, -16, 1)    &4  &  12&     64\ &  720 & 0 & 4 & 12& 40\\
69.  &\frac{1}{21}(13, -13, 1)    &8  &  21&     98\ &  1106 & 3 & 5 & 3& 4\\
70.  &\frac{1}{21}(11, -11, 1)    &2  &  6&     32\ &  399 & 0 & 2 & 6& 59\\
71.  &\frac{1}{22}(21, -21, 1)    &1  &  3&     16\ &  200 & 0 & 1 & 3& 30\\
72.  &\frac{1}{22}(19, -19, 1)    &7  &  21&     93\ &  1065 & 0 & 7 & 2& 8\\
73.  &\frac{1}{22}(17, -17, 1)    &9  &  22&     106\ &  1187 & 5 & 4 & 4& 3\\
74.  &\frac{1}{22}(15, -15, 1)    &3  &  9&     48\ &  574 & 0 & 3 & 9& 64\\
75.  &\frac{1}{22}(13, -13, 1)    &5  &  15&     77\ &  859 & 0 & 5 & 12& 30\\
76.  &\frac{1}{23}(22, -22, 1)    &1  &  3&     16\ &  200 & 0 & 1 & 3& 30\\
77.  &\frac{1}{23}(21, -21, 1)    &11  &  23&     114\ &  1282 & 10 & 1 & 1& 6\\
78.  &\frac{1}{23}(20, -20, 1)    &8  &  23&     102\ &  1165 & 1 & 7 & 1& 3\\
79.  &\frac{1}{23}(19, -19, 1)    &6  &  18&     88\ &  988 & 0 & 6 & 10& 24\\
80.  &\frac{1}{23}(18, -18, 1)    &9  &  23&     109\ &  1223 & 4 & 5 & 4& 2\\
81.  &\frac{1}{23}(17, -17, 1)    &4  &  12&     64\ &  734 & 0 & 4 & 12& 54\\
82.  &\frac{1}{23}(16, -16, 1)    &10  &  23&     112\ &  1263 & 7 & 3 & 3& 11\\
83.  &\frac{1}{23}(15, -15, 1)    &3  &  9&     48\ &  578 & 0 & 3 & 9& 68\\
84.  &\frac{1}{23}(14, -14, 1)    &5  &  15&     78\ & 870 & 0 & 5 &
13& 34\\\hline
\end{array}
$$
$$
\begin{array}{cccccccccc}\hline
No. \hspace*{2mm}&\hspace*{2mm}Singularity\hspace*{2mm} &\hspace*{2mm}\nabla_{1}\hspace*{2mm} &\hspace*{2mm}\nabla_{2}\hspace*{2mm} &\hspace*{2mm}\nabla_{3}\hspace*{2mm} \ &\hspace*{2mm}\nabla_{4}\hspace*{2mm} &\hspace*{2mm}\nabla'_{1} \hspace*{2mm}&\hspace*{2mm}\nabla'_{2} \hspace*{2mm}&\nabla'_{3} \hspace*{2mm}&\hspace*{2mm}\nabla'_{4}\\
85.  &\frac{1}{23}(13, -13, 1)    &7  &  21&     95\ &  1087 & 0 & 7 & 4& 16\\
86.  &\frac{1}{23}(12, -12, 1)    &2  &  6&     32\ &  400 & 0 & 2 & 6& 60\\
87.  &\frac{1}{24}(23, -23, 1)    &1  &  3&     16\ &  200 & 0 & 1 & 3& 30\\
88.  &\frac{1}{24}(19, -19, 1)    &5  &  15&     79\ &  880 & 0 & 5 & 14& 37\\
89.  &\frac{1}{24}(17, -17, 1)    &7  &  21&     97\ &  1107 & 0 & 7 & 6& 22\\
90.  &\frac{1}{24}(13, -13, 1)    &11  &  24&     118\ &  1332 & 9 & 2 & 2& 12\\
91.  &\frac{1}{25}(24, -24, 1)    &1  &  3&     16\ &  200 & 0 & 1 & 3& 30\\
92.  &\frac{1}{25}(23, -23, 1)    &12  &  25&     124\ &  1394 & 11 & 1 & 1& 6\\
93.  &\frac{1}{25}(22, -22, 1)    &8  &  24&     106\ &  1214 & 0 & 8 & 2& 8\\
94.  &\frac{1}{25}(21, -21, 1)    &6  &  18&     91\ &  1014 & 0 & 6 & 13& 29\\
95.  &\frac{1}{25}(19, -19, 1)    &4  &  12&     64\ &  745 & 0 & 4 & 12& 65\\
96.  &\frac{1}{25}(18, -18, 1)    &7  &  21&     99\ &  1125 & 0 & 7 & 8& 26\\
97.  &\frac{1}{25}(17, -17, 1)    &3  &  9&     48\ &  585 & 0 & 3 & 9& 75\\
98.  &\frac{1}{25}(16, -16, 1)    &11  &  25&     122\ &  1377 & 8 & 3 & 3& 13\\
99.  &\frac{1}{25}(14, -14, 1)    &9  &  25&     113\ &  1287 & 2 & 7 & 2& 6\\
100.  &\frac{1}{25}(13, -13, 1)    &2  &  6&     32\ &  400 & 0 & 2 & 6& 60\\
101. &\frac{1}{26}(25, -25, 1)& 1 &  3  &  16 \ &  200   &  0    &    1    &    3   &30\\
102. &\frac{1}{26}(23, -23, 1)    & 9    & 26  & 115 \ & 1314   &    1
&   8   &   1  &
3\\
103. &\frac{1}{26}(21, -21, 1)   & 5    & 15  &  80 \ &  898   &    0
&   5   &  15  &
48\\
104. &\frac{1}{26}(19, -19, 1)   &11    & 26  & 126 \ & 1417   &    7   &   4   &   4  &   9\\
105.  &\frac{1}{26}(17, -17, 1)   & 3    &  9  &  48 \ &  588   &
0 &   3   &   9  & 78\\\hline
\end{array}
$$

$$
\begin{array}{cccccccccc}\hline
No. \hspace*{2mm}&\hspace*{2mm}Singularity\hspace*{2mm} &\hspace*{2mm}\nabla_{1}\hspace*{2mm} &\hspace*{2mm}\nabla_{2}\hspace*{2mm} &\hspace*{2mm}\nabla_{3}\hspace*{2mm} \ &\hspace*{2mm}\nabla_{4}\hspace*{2mm} &\hspace*{2mm}\nabla'_{1} \hspace*{2mm}&\hspace*{2mm}\nabla'_{2} \hspace*{2mm}&\nabla'_{3} \hspace*{2mm}&\hspace*{2mm}\nabla'_{4}\\
106. &\frac{1}{26}(15 ,-15, 1)   & 7    & 21  & 101 \ & 1141   &    0
&   7   &  10  &
28\\
107. &\frac{1}{27}(26, -26, 1)    & 1    &  3  &  16 \ &  200   &    0
&   1   &   3  &
30\\
108.&\frac{1}{27}(25,-25, 1)     &13    & 27  & 134 \ & 1506   &   12
&   1   &   1  &
6\\
109. &\frac{1}{27}(23,-23, 1)    & 7    & 21  & 103\ & 1155    &    0
&   7   &  12  &
28\\
110. &\frac{1}{27}(22,-22, 1)    &11    & 27  & 130 \ & 1455      &
6     &  5     &   5    &
3\\
111.&\frac{1}{27}(20,-20, 1)     & 4    & 12  &  64 \ &  755     &
0     &       4     &  12    &
75\\
112.&\frac{1}{27}(19, -19, 1)      &10    & 27  & 124 \ & 1407      &
3     &       7     &   3    &
7\\
113.&\frac{1}{27}(17, -17, 1)      & 8    & 24  & 110 \ & 1257 &     0
&   8     &   6    &
23\\
114. &\frac{1}{27}(16, -16, 1)  & 5    & 15  &  80 \ &  906  &     0
& 5 &  15    &
56\\
115. &\frac{1}{27}(14, -14, 1)  & 2    &  6  &  32 \ &  400  &     0
& 2 &   6    &
60\\
116.  &\frac{1}{28}(27, -27, 1)   & 1    &  3  &  16 \ &  200  &     0
& 1 &   3    &
30\\
117.  &\frac{1}{28}(25, -25, 1)   & 9    & 27  &  119\ &  1363 &     0
& 9 &   2    &
8\\
118.  &\frac{1}{28}(23, -23, 1)   &11    & 28  &  133\ &  1491 &     5
& 6 &   5    &
2\\
119.   &\frac{1}{28}(19, -19, 1)  & 3    &   9 &   48\ &   593 &     0
&  3&    9   &
83\\
120   &\frac{1}{28}(17, -17, 1)  &5    &  15 &   80\ &   913   &
0    &    5    &   15   &
63\\
121.  &\frac{1}{28}(15, -15, 1)   &13    &  28 &   138\ &    1556    &
11    &     2    &    2   &
12\\
122. &\frac{1}{29}(28, -28, 1)    & 1    &   3 &    16\ &    200     &   0    &   1    &    3   &    30\\
123.  &\frac{1}{29}(27, -27, 1)   &14    &  29 &   144\ &    1618    &
13    &   1    &    1   &
6\\
124.  &\frac{1}{29}(26, -26, 1)   &10    &  29 &   128\ &  1463      &
1    &   9    &    1   &
3\\
125.  &\frac{1}{29}(25, -25, 1)   & 7    &  21 &   106\ &   1181     &
0    &   7    &   15   &
33\\
126.  &\frac{1}{29}(24, -24, 1)    &6    &  18 &    95\ &   1058
& 0    &   6    &   17   & 45\\\hline
\end{array}
$$
$$
\begin{array}{cccccccccc}\hline
No. \hspace*{2mm}&\hspace*{2mm}Singularity\hspace*{2mm} &\hspace*{2mm}\nabla_{1}\hspace*{2mm} &\hspace*{2mm}\nabla_{2}\hspace*{2mm} &\hspace*{2mm}\nabla_{3}\hspace*{2mm} \ &\hspace*{2mm}\nabla_{4}\hspace*{2mm} &\hspace*{2mm}\nabla'_{1} \hspace*{2mm}&\hspace*{2mm}\nabla'_{2} \hspace*{2mm}&\nabla'_{3} \hspace*{2mm}&\hspace*{2mm}\nabla'_{4}\\
127.   &\frac{1}{29}(23, 6, 1)   &5    &  15 &    80\ &    919     &
0    &   5    &   15   &
69\\
128.  &\frac{1}{29}(22, -22, 1)   &4    &  12 &    64\ &    764    &
0    &   4    &   12   &
84\\
129.  &\frac{1}{29}(21, -21, 1)   &11    &  29 &   135\ &   1525    &
4    &     7    &    4   &
6\\
130. &\frac{1}{29}(20, -20, 1)      &13    &  29&    142\ &   1604    &
10    & 3    &    3   &
16\\
131.  &\frac{1}{29}(19, -19, 1)     &3    &   9&     48\ &    595    &
0    & 3    &    9   &
85\\
132.  &\frac{1}{29}(18, -18, 1)     &8    &  24&    114\ &   1293    &
0    & 8    &   10   &
31\\
133.  &\frac{1}{29}(17, -17, 1)    &12    &  29  &  140\ &   1570    &
7    & 5    &    5   &
6\\
134.  &\frac{1}{29}(16, -16, 1)    & 9    &  27  &  121\ &   1385    &
0    & 9    &    4   &
16\\
135.   &\frac{1}{29}(15, -15, 1)   & 2    & 6    & 32  \ &  400      &
0    & 2    &    6   &
60\\
136.  &\frac{1}{30}(29, -29, 1)    & 1    & 3    & 16  \ &  200      &
0    & 1    &    3   &
30\\
137.  &\frac{1}{30}(23,-23, 1)    &13    &30    & 146 \ &  1646     &
9    &    4    &    4   &
14\\
138. &\frac{1}{30}(19, -19, 1)     &11    &30    & 137 \ &  1557     &
3    &    8    &    3   &
8\\
139.  &\frac{1}{30}(17, -17, 1)    & 7    &21    & 107 \ &  1193     &
0    &    7    &   16   &
38\\
140.  &\frac{1}{31}(30, -30, 1)     & 1    & 3    & 16  \ &   200     &
0    &    1    &    3   &
30\\
141.  &\frac{1}{31}(29, -29, 1)     &15   & 31   & 154 \ &   1730  &
14    &    1    &    1   &
6\\
142.  &\frac{1}{31}(28, -28, 1)     &10   & 30   & 132\ &   1512   &
0    &   10    &    2   &
8\\
143.   &\frac{1}{31}(27, -27, 1)    &8    &24    & 118\ &   1322   &
0    &   8    &   14   &
32\\
144.  &\frac{1}{31}(26, -26, 1)     &6    &18    & 96 \ &  1076    &
0    &   6    &   18   &
56\\
145.  &\frac{1}{31}(25, -25, 1)    &5    &15    & 80 \ &   930    &
0    &   5    &   15   &
80\\
146. &\frac{1}{31}(24, -24, 1)      &9    &27    &125 \ &  1426    &
0   &    9   &     8  &
29\\
147.  &\frac{1}{31}(23, -23, 1)     &4    &12    & 64 \ &   772    &
0   &    4   &    12  & 92\\\hline
\end{array}
$$
$$
\begin{array}{cccccccccc}\hline
No. \hspace*{2mm}&\hspace*{2mm}Singularity\hspace*{2mm} &\hspace*{2mm}\nabla_{1}\hspace*{2mm} &\hspace*{2mm}\nabla_{2}\hspace*{2mm} &\hspace*{2mm}\nabla_{3}\hspace*{2mm} \ &\hspace*{2mm}\nabla_{4}\hspace*{2mm} &\hspace*{2mm}\nabla'_{1} \hspace*{2mm}&\hspace*{2mm}\nabla'_{2} \hspace*{2mm}&\nabla'_{3} \hspace*{2mm}&\hspace*{2mm}\nabla'_{4}\\
148.  &\frac{1}{31}(22, -22, 1)      &7    &21    &108 \ &  1205    &
0   &    7   &    17  &
43\\
149.  &\frac{1}{31}(21, -21, 1)      &3    & 9    & 48 \ &   598
&        0   &    3   &     9  &
88\\
150.  &\frac{1}{31}(20, -20, 1)       &14   & 31   & 152\ &   1717
&       11   &    3   &     3  &
17\\
151.   &\frac{1}{31}(19, -19, 1)     &13   & 31   & 150\ &   1685
&        8   &    5   &     5  &
9\\
152.  &\frac{1}{31}(18, -18, 1)      &12   & 31   & 146\ &   1642
&        5   &    7   &     5  &
4\\
153.  &\frac{1}{31}(17, -17, 1)      &11   & 31   & 139\ &   1585
&        2   &     9   &     2  &
6\\
154. &\frac{1}{31}(16, -16, 1)       &2    & 6    &32  \ &  400
&        0   &     2   &     6  &
60\\
155.  &\frac{1}{32}(31, -31, 1)      &1    & 3    &16  \ &  200    &
0    &    1    &    3   &
30\\
156. &\frac{1}{32}(29, -29, 1)      &11    &32    &141 \  & 1612   &
1    &   10    &    1   &
3\\
157. &\frac{1}{32}(27, -27, 1)      &13    &32    &154 \ &   1723     &
7    &    6    &    6   &
3\\
158. &\frac{1}{32}(25, -25, 1)      & 9    &27    &127 \ & 1444 & 0 & 9
&   10   &
33\\
159.  &\frac{1}{32}(23, -23, 1)     & 7    &21    &109\ &   1216      &
0    &    7    &   18   &
47\\
160. &\frac{1}{32}(21, -21, 1)      & 3    & 9    & 48\ &    599
&     0    &    3    &    9   &
89\\
161. &\frac{1}{32}(19, -19, 1)      & 5    &15    & 80\ &    935
&     0    &    5    &   15   &
85\\
162.&\frac{1}{32}(17, -17, 1)       &15    &32    &158\ &  1780     &
13    &    2    &    2   &
12\\
163. &\frac{1}{33}(32, -32, 1)     & 1    & 3    &16 \ &  200      &
0    &    1    &    3   &
30\\
164.  &\frac{1}{33}(31, -31, 1)     &16    &33    &164\ &  1842     &
15    &    1    &    1   &
6\\
165.  &\frac{1}{33}(29, -29, 1)     & 8    &24    &121\ &  1348     &
0    &    8    &   17   &
37\\
166.  &\frac{1}{33}(28, -28, 1)     &13    &33    & 157\ &   1759   &
6    &    7    &    6   &
2\\
167.   &\frac{1}{33}(26, -26, 1)    &14    &33    &160\ &  1800     &
9    &    5    &    5   &
12\\
168.  &\frac{1}{33}(25, -25, 1)     & 4    &12    &64 \ &  779
& 0    &    4    &   12   & 99\\\hline
\end{array}
$$
$$
\begin{array}{cccccccccc}\hline
No. \hspace*{2mm}&\hspace*{2mm}Singularity\hspace*{2mm} &\hspace*{2mm}\nabla_{1}\hspace*{2mm} &\hspace*{2mm}\nabla_{2}\hspace*{2mm} &\hspace*{2mm}\nabla_{3}\hspace*{2mm} \ &\hspace*{2mm}\nabla_{4}\hspace*{2mm} &\hspace*{2mm}\nabla'_{1} \hspace*{2mm}&\hspace*{2mm}\nabla'_{2} \hspace*{2mm}&\nabla'_{3} \hspace*{2mm}&\hspace*{2mm}\nabla'_{4}\\
169.  &\frac{1}{33}(23, -23, 1)     &10    &30   &136 \ & 1556      &
0    &   10    &    6   &
24\\
170. &\frac{1}{33}(20, -20, 1)       & 5    &15    &80 \ &  940     &
0    &    5    &   15   &
90\\
171.  &\frac{1}{33}(19, -19, 1)     & 7    &21  & 110 \ &   1226   &
0    &    7    &   19   &
50\\
172.  &\frac{1}{33}(17, -17, 1)     & 2     &6  &  32 \ &    400   &
0    &    2    &    6   &
60\\
173.  &\frac{1}{34}(33, -33, 1)     & 1     &3  &  16 \ &    200   &
0    &    1    &    3   &
30\\
174.  &\frac{1}{34}(31, -31, 1)     &11  & 33   &  145\ &    1661  &
0    &   11    &    2   &
8\\
175.   &\frac{1}{34}(29, -29, 1)    & 7   &21   &  111\ &    1236     &
0    &    7    &   20   &
53\\
176.  &\frac{1}{34}(27, -27, 1)     & 5   &15   &  80 \ &    945      &
0    &    5    &   15   &
95\\
177.  &\frac{1}{34}(25, -25, 1)     &15   &34   &   166\ &    1874    &
11    &    4    &    4   &
18\\
178. &\frac{1}{34}(23, -23, 1)      & 3   & 9   & 48  \ &  600       &
0    &    3    &    9   &
90\\
179.  &\frac{1}{34}(21, -21, 1)     &13   &34   & 159  \ &  1793     &
5    &    8    &    5   &
6\\
180.  &\frac{1}{34}(19, -19, 1)     & 9   &27   & 131  \ &  1475     &
0    &    9    &   14   &
36\\
181.  &\frac{1}{35}(34, -34, 1)     & 1   & 3   & 16   \ &  200      &
0    &    1    &    3   &
30\\
182.  &\frac{1}{35}(33,-33, 1)     &17   &35   & 174  \ &  1954     &
16     &   1     &   1    &
6\\
183.   &\frac{1}{35}(32, -32, 1)    &12   &35   & 154  \ &  1761   &
1     &  11     &   1    &
3\\
184.  &\frac{1}{35}(31, -31, 1)     & 9   &27   & 133  \ &     1489&
0     &   9     &  16    &
36\\
185.  &\frac{1}{35}(29, -29, 1)     & 6   &18   & 96   \ &1104    &   0
&   6     &  18    &
84\\
186. &\frac{1}{35}(27, -27, 1)     &13   &35   & 161  \ &1826    &   4
&   9     &   4    &
9\\
187.  &\frac{1}{35}(26, -26, 1)     & 4   &12   & 64   \ &785     & 0 &
4     &  12    &
105\\
188.  &\frac{1}{35}(24, -24, 1)    &16   &35   & 172  \ &1942    &
13     &   3     &   3    &
18\\
189.  &\frac{1}{35}(23, -23, 1)    & 3   & 9   & 48   \ &600     &
0 &   3     &   9    & 90\\\hline
\end{array}
$$
$$
\begin{array}{cccccccccc}\hline
No. \hspace*{2mm}&\hspace*{2mm}Singularity\hspace*{2mm} &\hspace*{2mm}\nabla_{1}\hspace*{2mm} &\hspace*{2mm}\nabla_{2}\hspace*{2mm} &\hspace*{2mm}\nabla_{3}\hspace*{2mm} \ &\hspace*{2mm}\nabla_{4}\hspace*{2mm} &\hspace*{2mm}\nabla'_{1} \hspace*{2mm}&\hspace*{2mm}\nabla'_{2} \hspace*{2mm}&\nabla'_{3} \hspace*{2mm}&\hspace*{2mm}\nabla'_{4}\\
190.  &\frac{1}{35}(22, -22, 1)    & 8   &24   & 123  \ &1372    &   0
&   8     &  19    &
47\\
191.   &\frac{1}{35}(19, -19, 1)   &11   &33  &  147  \ &1683   &    0
&   11    &    4   &
16\\
192.  &\frac{1}{35}(18, -18, 1)   & 2   & 6  &  32   \ &400    &    0
&    2    &    6   &
60\\
193.  &\frac{1}{36}(35, -35, 1)    & 1   & 3  &   16  \ & 200   &    0
&    1    &    3   &
30\\
194. &\frac{1}{36}(31, -31, 1)     & 7  & 21   & 112  \ & 1254  &    0
&    7    &   21   &
64\\
195.  &\frac{1}{36}(29, -29, 1)    & 5  & 15  &  80   \ &954    &    0
&    5    &   15   &
104\\
196.  &\frac{1}{36}(25, -25, 1)    &13  &  36 &   163 \ & 1856  &    3
&   10    &    3   &
9\\
197.  &\frac{1}{36}(23, -23, 1)    &11  &  33 &   149 \ & 1705  &    0
&   11    &    6   &
24\\
198.  &\frac{1}{36}(19, -19, 1)    &17  &  36 &   178 \ & 2004  &   15
&    2    &    2   &
12\\
199.   &\frac{1}{37}(36, -36, 1)   & 1  &   3 &   16  \ & 200      & 0
&    1    &    3   &
30\\
200. &\frac{1}{37}(35, -35, 1)    &18  &  37 &   184 \ & 2066     &17
&    1    &    1   &
6\\
201. &\frac{1}{37}(34, -34, 1)    &12  &  36 &   158 \ & 1810     & 0
&   12    &    2   &
8\\
202.&\frac{1}{37}(33, -33, 1)     & 9  &  27 &   136 \ & 1515     & 0
&    9    &   19   &
41\\
203. &\frac{1}{37}(32, -32, 1)    &15  &  37 &   178 \ & 1991     & 8
&    7    &    7   &
3\\
204.  &\frac{1}{37}(31, -31, 1)   & 6  &  18 &   96  \ &1115    &
0    &    6    &   18   &
95\\
205.  &\frac{1}{37}(30, -30, 1)   &16  &  37 &   180 \ & 2029   &
11    &    5    &    5   &
17\\
206.  &\frac{1}{37}(29, -29, 1)   &14  &  37&    172 \ & 1944   &
5    &    9    &    5   &
8\\
207.   &\frac{1}{37}(28, -28, 1)  & 4  &  12&     64 \ &  790   &
0    &    4    &   12   &
110\\
208.  &\frac{1}{37}(27, -27, 1)   &11  &  33&    151 \ & 1726   &
0    &   11    &    8   &
31\\
209.  &\frac{1}{37}(26, -26, 1)     &10  &  30&    144 \ & 1628
&   0    &  10    &   14   &
40\\
210. &\frac{1}{37}(25, -25, 1)      & 3  &   9&     48 \ &  600 &
0    &   3    &    9   & 90\\\hline
\end{array}
$$
$$
\begin{array}{cccccccccc}\hline
No. \hspace*{2mm}&\hspace*{2mm}Singularity\hspace*{2mm} &\hspace*{2mm}\nabla_{1}\hspace*{2mm} &\hspace*{2mm}\nabla_{2}\hspace*{2mm} &\hspace*{2mm}\nabla_{3}\hspace*{2mm} \ &\hspace*{2mm}\nabla_{4}\hspace*{2mm} &\hspace*{2mm}\nabla'_{1} \hspace*{2mm}&\hspace*{2mm}\nabla'_{2} \hspace*{2mm}&\nabla'_{3} \hspace*{2mm}&\hspace*{2mm}\nabla'_{4}\\
211.  &\frac{1}{37}(24, -24, 1)     &17  &  37&     182\ &  2054
&  14    &   3    &    3   &
18\\
212.  &\frac{1}{37}(23, -23, 1)     & 8  &  24&      125\ &   1394  &
0    &   8    &   21   &
55\\
213.  &\frac{1}{37}(22, -22, 1)     & 5  &  15&      80\ &   958    &
0    &   5    &   15   &
108\\
214.  &\frac{1}{37}(21, -21, 1)     & 7  &  21&     112\ &  1262    &
0    &   7    &   21   &
72\\
215.   &\frac{1}{37}(20, -20, 1)    &13  &  37&     165\ &  1883    &
2    &  11    &    2   &
6\\
216.  &\frac{1}{37}(19, -19, 1)     & 2  &   6&      32\ &   400 & 0
&   2    &    6   &
60\\
217.  &\frac{1}{38}(37, -37, 1)     & 1  &   3&      16\ &   200
&        0    &   1    &    3   &
30\\
218. &\frac{1}{38}(35, -35, 1)    &13  &  38&     167\ &  1910      &
1    &       12    &    1   &
3\\
219.  &\frac{1}{38}(33, -33, 1)     &15  &  38&     181\ & 2027
&        7    &        8    &    7   &
2\\
220.  &\frac{1}{38}(31, -31, 1)   &11  &  33&     153\ &   1745    &
0    &       11    &   10   &
36\\
221.  &\frac{1}{38}(29, -29, 1)    &17  &  38&     186\ &  2101     &
13    &   4    &    4   &
21\\
222.  &\frac{1}{38}(27, -27, 1)    & 7  &  21&     112\ &  1270     &
0    &   7    &   21   &
80\\
223.   &\frac{1}{38}(25, -25, 1)  & 3  &   9&      48\ &   600     &
0    &   3    &    9   &
90\\
224.  &\frac{1}{38}(23, -23, 1)    & 5  &  15&      80\ &   962     &
0    &   5    &   15   &
112\\
225.  &\frac{1}{38}(21, -21, 1)    & 9  &  27&     137\ &  1527   &
0    &   9    &   20   &
46\\
226. &\frac{1}{39}(38, -38, 1)     & 1  &   3&      16\ &   200   &
0    &   1    &    3   &
30\\
227.  &\frac{1}{39}(37, 2, 1)    &19  &  39&     194\ &  2178   &
18    &   1    &    1   &
6\\
228.  &\frac{1}{39}(35, -35, 1)    &10  &  30&     148\ &  1656   &
0    &  10    &   18   &
40\\
229.  &\frac{1}{39}(34, -34, 1)     & 8  &  24&     127\ &  1414   &
0    &   8    &   23   &
61\\
230.  &\frac{1}{39}(32, -32, 1)     &11  &  33&     155\ &  1763
& 0    &  11    &   12   &
40\\
231.   &\frac{1}{39}(31, -31, 1)    & 5  &  15&      80\ &   966 & 0
&   5    &   15   & 116\\\hline
\end{array}
$$
$$
\begin{array}{cccccccccc}\hline
No. \hspace*{2mm}&\hspace*{2mm}Singularity\hspace*{2mm} &\hspace*{2mm}\nabla_{1}\hspace*{2mm} &\hspace*{2mm}\nabla_{2}\hspace*{2mm} &\hspace*{2mm}\nabla_{3}\hspace*{2mm} \ &\hspace*{2mm}\nabla_{4}\hspace*{2mm} &\hspace*{2mm}\nabla'_{1} \hspace*{2mm}&\hspace*{2mm}\nabla'_{2} \hspace*{2mm}&\nabla'_{3} \hspace*{2mm}&\hspace*{2mm}\nabla'_{4}\\
232.  &\frac{1}{39}(29, -29, 1)     & 4  &  12&       64\ &    794
& 0    &   4    &   12   &
114\\
233.  &\frac{1}{39}(28, 11, 1)     & 7  &  21&    112 \ & 1277  &
0    &    7    &   21   &
87\\
234. &\frac{1}{39}(25, -25, 1)      &14  &  39&    176 \ & 2005  &
3    &   11    &    3   &
9\\
235.  &\frac{1}{39}(23, -23, 1)     &17  &  39&    190 \ & 2143  &
12    &    5    &    5   &
19\\
236.   &\frac{1}{39}(22, -22, 1)    &16  &  39  &  188  \ &2106   & 9
&    7    &    7   &
6\\
237.   &\frac{1}{39}(20, -20, 1)    & 2  &   6&     32 \ &  400  &
0    &    2    &    6   &
60\\
238.   &\frac{1}{40}(39, -39, 1)    & 1  &   3&     16 \ &  200  &
0    &    1    &    3   &
30\\
239.    &\frac{1}{40}(37, -37, 1)   &13  &  39&     171\ &  1959
&        0    &   13    &    2   &
8\\
240.   &\frac{1}{40}(33, -33, 1)    &17  &  40&     194\ &  2183   &
11    &    6    &    6   &
15\\
241.   &\frac{1}{40}(31, -31, 1)    & 9  &  27&    139 \ & 1551      &
0    &    9    &   22   &
56\\
242.  &\frac{1}{40}(29, -29, 1)    &11  &  33&    157 \ & 1780      &
0    &   11    &   14   &
43\\
243.   &\frac{1}{40}(27, -27, 1)    & 3  &   9&     48 \ &  600
&   0    &    3    &    9   &    90
\\
244.  &\frac{1}{40}(23, -23, 1)     & 7  &  21&     112\ &  1283
&   0    &    7    &   21   &
93\\
245.  &\frac{1}{40}(21, -21, 1)    &19  &  40&     198\ &  2228 & 17
&    2    &    2   &
12\\
246.  &\frac{1}{41}(40, -40, 1)     & 1  &   3&     16 \ &  200    &
0    &    1    &    3   &    30
\\
247.   &\frac{1}{41}(39, -39, 1)    &20  &  41&   204  \ &2290     &
19    &    1    &    1   &
6\\
248.  &\frac{1}{41}(38, -38, 1)     &14  &  41&    180 \ & 2059    &
1    &   13    &    1   &
3\\
249.  &\frac{1}{41}(37, -37, 1)     &10  &  30&    151 \ & 1682 &  0
&       10    &   21   &
45\\
250. &\frac{1}{41}(36, -36, 1)    & 8  &  24&    128 \ & 1432   &  0
&  8    &   24   &
72\\
251.  &\frac{1}{41}(35, -35, 1)    & 7  &  21&    112 \ & 1289   &  0
&  7    &   21   &
99\\
252. &\frac{1}{41}(34, -34, 1)      & 6  &  18&     96 \ & 1135   &
0 &  6    &   18   & 115\\\hline
\end{array}
$$
$$
\begin{array}{cccccccccc}
\hline
No. \hspace*{2mm}&\hspace*{2mm}Singularity\hspace*{2mm} &\hspace*{2mm}\nabla_{1}\hspace*{2mm} &\hspace*{2mm}\nabla_{2}\hspace*{2mm} &\hspace*{2mm}\nabla_{3}\hspace*{2mm} \ &\hspace*{2mm}\nabla_{4}\hspace*{2mm} &\hspace*{2mm}\nabla'_{1} \hspace*{2mm}&\hspace*{2mm}\nabla'_{2} \hspace*{2mm}&\nabla'_{3} \hspace*{2mm}&\hspace*{2mm}\nabla'_{4}\\
253.  &\frac{1}{41}(33, -33, 1)     & 5  &  15&     80 \ & 973 & 0 &
5    &   15   &
123\\
254.  &\frac{1}{41}(32, -32, 1)     & 9  &  27&    140 \ & 1562      &
0    &  9    &   23   &
60\\
255.  &\frac{1}{41}(31, -31, 1)     & 4  &  12&     64 \ &  797
&  0    &        4    &   12   &
117\\
256.  &\frac{1}{41}(30, -30, 1)    &15  &  41&    187 \ & 2126      &
4    &  11    &    4   &
11\\
257.   &\frac{1}{41}(29, -29, 1)   &17  &  41&    198 \ & 2221      &
10    &   7    &    7   &
9\\
258.  &\frac{1}{41}(28,-28, 1)    &19  &  41&    202 \ & 2278      &
16    &   3    &    3   &
18\\
259.  &\frac{1}{41}(27,-27, 1)    & 3  &   9&     48 \ &  600      &
0    &   3    &    9   &
90\\
260. &\frac{1}{41}(26, -26, 1)      &11  &  33&    159 \ & 1795
& 0    &  11    &   16   &
44\\
261.  &\frac{1}{41}(25, -25, 1)     &18  &  41&    200 \ & 2257  &
13    &   5    &    5   &
21\\
262.   &\frac{1}{41}(24, -24, 1)    &12  &  36&    166 \ & 1895     &
0    &  12    &   10   &
37\\
263.   &\frac{1}{41}(23, -23, 1)    &16  &  41&    194 \ & 2178     &
7    &   9    &    7   &
4\\
264.   &\frac{1}{41}(22, -22, 1)    &13  &  39&    173\ & 1981      &
0    &  13    &    4   &
16\\
265.    &\frac{1}{41}(21, -21, 1)   & 2  &   6&     32\ &   400     &
0    &   2    &    6   &
60\\
266.   &\frac{1}{42}(41,-41, 1)    & 1  &   3&     16\ &   200     &
0    &   1    &    3   &
30\\
267.   &\frac{1}{42}(37, -37, 1)    &17  &  42&    202\ & 2259    &
9    &   8    &    8   &
3\\
268.  &\frac{1}{42}(31, -31, 1)    &19  &  42   & 206 \ & 2327   &
5    &   4    &    4   &
23\\
269.   &\frac{1}{42}(29, -29, 1)    &13  &  39   &175  \ & 2003
&  0 &  13    &    6   &
24\\
270.  &\frac{1}{42}(25, -25, 1)     & 5  &  15  &  80 \ &  976     &
0    &        5    &   15   &
126\\
271.  &\frac{1}{42}(23,-23,1)     &11  &  33  & 161 \ & 1809 & 0 &
11    &   18   &    44 \\\hline
\end{array}
$$\vspace*{0.2cm}\hspace*{1cm}

$\star$The first 100 datum can also be found in \cite{2}. And to my
knowledge, the notations $\Delta_{n}$ were firstly used by Fletcher
and Reid (see \cite{2}).


\clearpage

\begin{thebibliography}{99}
\bibitem[1]{Bom} E. Bombieri, {\it Canonical models of surfaces of general type}. Publ. I.H.E.S. {\bf 42} (1973), 171--219.

\bibitem[2]{cat}F. Catanese, {\it Pluricanonical mappings of surfaces with $K^{2}=1, 2$,
$q=p_{g}=0$}. C.I.M.E. 1977: Algebraic Surfaces, Liguori,
Napoli(1981) 247-266.

\bibitem[3]{C-H}J.A. Chen and C.D. Hacon, {\it Linear series of irregular varieties}. Algebraic Geometry in
East Asia, Japan, 2002, World Scientific Press.

\bibitem[4]{0}J.A. Chen and M. Chen, {\it The Canonical volume of 3-folds of general type with $\chi\le 0$}. arXiv: 07041702.

\bibitem[5]{7}J.A. Chen and M. Chen, {\it Explicit birational geometry of threefolds of general
type}. (2007), arXiv:0706.2987v2 [math.AG].

\bibitem[6]{5}M. Chen, {\it Canonical stability in terms of singularity index for algebraic threefolds}. Math. Proc. Camb. Phil. Soc. {\bf 131}
(2001), 241-264.

\bibitem[7]{11}M. Chen, {\it A sharp lower bound for the canonical volume of 3-folds of general type}. Math. Annalen 337 (2007), {\bf no.4}, 887-908.

\bibitem[8]{0000}M. Chen and D.Q. Zhang, {\it Characterization of the 4-canonical birationality of algebraic threefolds}. Math. Z. DOI 10.1007/s00209-007-0186-4.

\bibitem[9]{chenzhu}M. Chen and L. Zhu, {\it Projective 3-folds of general type with
$\chi=1$}. Proceeding ICCM 2007(to appear). arXiv:math/0611898v2.

\bibitem[10]{8}C. Ciliberto, {\it The bicanonical map for surfaces of general
type}. Proc. Symposia in Pure Math. {\bf 62} (1997), 57-83.

\bibitem[11]{10}A. Corti and M. Reid, {\it Explicit birational geometry of 3-folds}. London Mathematical Society, Lecture Note Series, {\bf 281} (2000), Cambridge University Press,
Cambridge.

\bibitem[12]{2}A.R. Fletcher, {\it Contributions to Riemann-Roch on Projective 3-folds with Only Canonical
Singularities and Applications}. Proceedings of Symposia in Pure
Mathematics {\bf 46} (1987), 221-231.

\bibitem[13]{Fletcher} A.R. Iano-Fletcher,
{\it Inverting Reid's exact plurigenera formula}. Math. Ann. {\bf
284} (1989), no. 4, 617-629.

\bibitem[14]{H-M} C.D. Hacon and J. M$^{\text{\rm c}}$Kernan, {\it
Boundedness of pluricanonical maps of varieties of general type}.
Invent. Math. {\bf 166} (2006), 1-25.

\bibitem[15]{3}H. Hironaka, {\it Resolution of
singularities of an algebraic variety over a field of characteristic
zero I}. Ann. of Math. {\bf 79} (1964), 109-203; II, ibid. 205-326.

\bibitem[16]{6}J. Koll\'ar, {\it
Higher direct images of dualizing sheaves I}. Ann. of Math. {\bf
123} (1986), 11-42.

\bibitem[17]{Miyao}Y. Miyaoka, {\it Tricanonical maps of numerical Godeaux
surfaces}. Invent. Math. {\bf 34} (1976) 99-111.

\bibitem[18]{133}Y. Miyaoka, {\it The pseudo-effectivity of $3c_2-c_1^2$ for
varieties with numerically effective canonical classes}. Algebraic
Geometry, Sendai, 1985. Adv. Stud. Pure Math. {\bf 10}(1987),
449-476.

\bibitem[19]{12}M. Reid, {\it Young person's guide to canonical
singularities}. Proc. Sympos. Pure Math. {\bf 46} (1987), 345-414.

\bibitem[20]{reider}I. Reider, {\it Vector bundles of rank 2 and linear systems on algebraic
surfaces}. Ann. Math. {\bf 127} (1988), 309-316.

\bibitem[21]{shin}D.K. Shin, {\it On a computation of plurigenera of a canonical
threefold}. J. of Alg. {\bf 309}(2007), 559-568.

\bibitem[22]{Tak} S. Takayama, {\it Pluricanonical systems on algebraic
varieties of general type}. Invent. Math. {\bf 165} (2006), 551-584.
\end{thebibliography}
\end{document}